\documentclass{article}    
\usepackage{amssymb, latexsym, amsmath, amsthm}

\newcommand{\bi} {\begin{itemize}}
\newcommand{\ei} {\end{itemize}}
\newcommand{\bea} {\begin{eqnarray}}
\newcommand{\eea} {\end{eqnarray}}
\newcommand{\be} {\begin{equation}}
\newcommand{\ee} {\end{equation}}
\newcommand{\bean} {\begin{eqnarray*}}
\newcommand{\eean} {\end{eqnarray*}}

\RequirePackage{epsfig}
\RequirePackage{amssymb}
\RequirePackage{natbib}
\RequirePackage{graphics}
\bibpunct{(}{)}{,}{a}{,}{,}

\begin{document}
\title{Bounded Point Evaluations For Certain Polynomial And Rational Modules}
\author{Liming Yang}
\author{
 Liming Yang \\
 Department of Mathematics \\
 Virginia Polytechnic and State University \\
 Blacksburg, VA 24061 \\
yliming@vt.edu
     } 
						
\date{}				

\maketitle

\newtheorem{Theorem}{Theorem}
\newtheorem*{MTheorem}{Main Theorem}
\newtheorem{Corollary}{Corollary}
\newtheorem{Definition}{Definition}
\newtheorem{Assumption}{Assumption}
\newtheorem{Lemma}{Lemma}
\newtheorem{Problem}{Problem}
\newtheorem*{Example}{Example}
\newtheorem{KnownResult}{Known Result}
\newtheorem{Algorithm}{Algorithm}
\newtheorem{Property}{Property}
\newtheorem{Proposition}{Proposition}


\abstract{
Let $K$ be a compact subset of the complex plane $\mathbb C.$ Let $P(K)$
and $R(K)$ be the closures in $C(K)$ of analytic polynomials and rational functions with
poles off $K,$ respectively. Let $A(K) \subset C(K)$ be the algebra of functions that are analytic in the interior of $K$. For $1\le t <\infty,$ let $P^t(1, \phi_1,...,\phi_N,K)$ be the closure of $P(K)+P(K)\phi_1+...+P(K)\phi_N$ in $L^t(dA|_K),$ where $dA|_K$ is the area measure restricted to $K$ and $\phi_1,...,\phi_N\in L^t(dA|_K).$ Let $HP(\phi_1,...,\phi_N,K)$ be the closure of $P(K)\phi_1+...+P(K)\phi_N +R(K)$ in $C(K),$ where $\phi_1,...,\phi_N\in C(K).$
In this paper, we prove if $R(K)\ne C(K),$ then there exists an analytic bounded point evaluation for both $P^t(1, \phi_1,...,\phi_N,K)$ and $HP(\phi_1,...,\phi_N,K)$ for certain smooth functions $\phi_1,...,\phi_N,$ in particular, for $\bar z,\bar z^2,...,\bar z^N.$  
We show that $A(K)\subset HP(\bar z,\bar z^2,...,\bar z^N,K)$ if and only if $R(K) = A(K).$ In particular, $C(K) \ne HP(\bar z,\bar z^2,...,\bar z^N,K)$ unless $R(K) = C(K).$ We also give an example of $K$ showing the results are not valid if we replace $\bar z^n$ by certain $\phi_n,$ that is, there exist $K$ and a function $\phi\in A(K)$ such that $R(K) \ne A(K),$ but $A(K) = HP (\phi ,K).$

}

\section{Introduction}

Let $\mathcal{P}$ denote the set of polynomials in the complex variable $z.$ 
For a compact subset $K$ of the complex plane $\mathbb C,$ let $Rat(K)$ be the set of all rational functions with poles off $K$ and let $C(K)$ denote the Banach algebra of complex-valued continuous functions on $K$ with customary norm $\|.\|_K.$ Let $P(K)$ and $R(K)$ denote the closures in $C(K)$ of $\mathcal{P}$ and $Rat(K),$ respectively. Let $A(K) \subset C(K)$ be the algebra of functions that are analytic in the interior of $K.$
For $\phi_1,...,\phi_N\in C(K),$ let $HP(\phi_1,...,\phi_N,K)$ denote the closure of $P(K)\phi_1+...+P(K)\phi_N +R(K)$ in $C(K).$
For $1\le t <\infty,$ let $L^t(K) = L^t(dA|_K),$ where $dA|_K$ is the area measure restricted to $K.$ For $\phi_1,...,\phi_N\in L^t(K),$ let $P^t(1, \phi_1,...,\phi_N,K)$ be the closure of $P(K)+P(K)\phi_1+...+P(K)\phi_N$ in $L^t(K).$
For a subset $A \subset \mathbb C,$ we set $Int(A)$ for its interior, $\bar A$ or $clos (A)$  for its closure, $A^c$ for its complement, and $\chi _A$ for its characteristic function.

For a subspace $A$ of $C(K)$ and a function $f\in C(K),$ we define the distance from $f$ to $A$ by
 \[
 \ dist(f, A) = \inf_{g\in A} \| f - g \| _ K.
 \]
For a subspace $B$ of $L^t(K)$ and a function $f\in L^t(K),$ we define the distance from $f$ to $B$ by
 \[
 \ dist(f, B) = \inf_{g\in B} \| f - g \| _ {L^t(K)}.
 \]
Set
 \[
 \ B (\lambda _0, \delta ) = \{z: |z - \lambda_0| < \delta \}.
 \]
The open unit disk is denoted by $\mathbb D =B (0, 1). $ The constants used in the paper such as $C,$ $C_0,$ $C_1,$ $C_N,$ $\delta_0,$ $\delta_1,$ $\delta_N,$ $\epsilon_0,$ $\epsilon_1,$ $\epsilon_N,...$ may change from one step to the next.

We denote the Riemann sphere $\mathbb C_\infty = \mathbb C \cup \{\infty \}.$ For a compact subset $E \subset \mathbb C,$ we
define the analytic capacity of $E$ by
\[
\ \gamma(E) = \sup_{f\in \mathcal{A} (E)} |f'(\infty)|,
\]
where $\mathcal{A} (E)$ consists of those functions $f$ analytic in $\mathbb C_\infty \setminus E$ for which $f(\infty)=0,$
$|f(z)| \le 1$ for all $z \in \mathbb C_\infty \setminus E,$ and
\[
\ f'(\infty) = \lim _{z \rightarrow \infty} z(f(z) - f(\infty)).
\]
The analytic capacity of a general $E_1 \subset \mathbb C$ is defined to be 
\[
\ \gamma (E_1) = \sup \{\gamma (E) : E \subset E_1, ~E~ compact\}.
\]
The continuous analytic capacity for a compact subset $E$ is defined similarly as 
 \[
\ \alpha(E) = \sup_{f\in \mathcal{AC} (E)} |f'(\infty)|,
\]
where $\mathcal{AC} (E) = \mathcal{A} (E) \cap C(\mathbb C_{\infty}).$ For a general $E_1 \subset \mathbb C,$
\[
\ \alpha (E_1) = \sup \{\alpha (E) : E \subset E_1, ~E~ compact\}.
\]
(see \cite{gamelin} and \cite{conway} for basic information of rational approximation and analytic capacity).

Let $\nu$ be a compactly supported finite measure on $\mathbb {C}.$ The Cauchy transform
of $\nu$ is defined by
\[
\ \mathcal C\nu (z) = \int \dfrac{1}{w - z} d\nu (w)
\]
for all $z\in\mathbb {C}$ for which
$\int \frac{d|\nu|(w)}{|w-z|} < \infty .$ A standard application of Fubini's
Theorem shows that $\mathcal C\nu \in L^s_{loc}(\mathbb {C} )$ for $ 0 < s < 2,$ in particular that it is
defined for Area almost all $z,$ and clearly $\mathcal C\nu$ is analytic in $\mathbb C_\infty \setminus spt \nu .$

We denote the map
 \[
 \ E^t(\lambda ) : p_0 + \sum_{i=1}^N p_i\phi_i \rightarrow \left[ \begin{array}{c}
p_0(\lambda )\\
p_1(\lambda ) \\
... \\
p_N(\lambda ) \end{array} \right],\tag{1-1}
 \]
where $p_0,p_1,...,p_N\in \mathcal P.$ If $E(\lambda)$ is bounded from $P^t(1, \phi_1,...,\phi_N,K)$ to $(\mathbb C^{N+1}, \|.\|_{N+1}),$ where $\|x\|_{N+1} = \sum_{i=0}^N |x_i|$ for $x\in \mathbb C^{N+1},$
 then every component in the right hand side extends to a bounded linear functional on $P^t(1, \phi_1,...,\phi_N,K)$ and we will call $\lambda$ a
bounded point evaluation for $P^t(1, \phi_1,...,\phi_N,K).$ 
A bounded point evaluation $\lambda_0$ is called an analytic bounded point evaluation for $P^t(1, \phi_1,...,\phi_N,K)$ if there is a neighborhood $B(\lambda_0, \delta)$ of $\lambda_0$ such that every $\lambda \in B(\lambda_0, \delta)$ is a bounded point evaluation and $E^t(\lambda)$ is analytic as a function of $\lambda$ on $B(\lambda_0, \delta)$ (equivalently (1-1) is uniformly bounded for  $\lambda\in B(\lambda_0, \delta)$).  Similarly, we can define a bounded point evaluation (or analytic bounded point evaluation) $\lambda$ for $HP(\phi_1,...,\phi_N,K)$by replacing (1-1) with the following map:
 \[
 \ E(\lambda ) : r + \sum_{i=1}^N p_i\phi_i \rightarrow \left[ \begin{array}{c}
p_1(\lambda )\\
p_2(\lambda ) \\
... \\
p_N(\lambda ) \end{array} \right],\tag{1-2}
 \]
where $p_1,p_2,...,p_N\in \mathcal P$ and $r\in Rat(K).$ Notice that the rational function $r\in Rat(K)$ does not appear on the right hand side of definition (1-2).

For an arbitrary finite compactly supported positive measure $\mu ,$ \cite{thomson} describes completely the structure of $P^t(\mu ),$ the closed
subspace of $L^t(\mu )$ spanned by $\mathcal P.$ \cite{ce93} extends some results of Thomson's Theorem to the space $R^t(K,\mu ),$ the closure of $Rat(K)$ in $L^t(\mu ),$ while \cite{b08} expresses $R^t(K,\mu )$ as a direct sum that includes both Thomson's theorem and results of \cite{ce93}. For a compactly supported complex Borel measure $\nu$ of $\mathbb C,$ by estimating analytic capacity of the set $\{\lambda: |\mathcal C\nu (\lambda)| \ge c \},$  \cite{b06}, \cite{ars}, and \cite{ARS10} provide interesting alternative proofs of Thomson's theorem. Both their proofs rely on X. Tolsa's deep results on analytic capacity. \cite{y17} extends some results to a rationally multicyclic subnormal operator (restriction of a normal operator on a separable Hilbert space to an invariant subspace).

However, even for $\mu = dA |_K,$ it is difficult to obtain necessary and sufficient conditions under which $P^t(1,K) = L^t(K)$ or $P^t(1,K)$ has a bounded point evaluation. \cite{bm} proves if $R(K) \ne C(K),$ then $P^t(1,K)$ has a bounded point evaluation. \cite{yang3} shows that there exists a compact subset $K\subset \mathbb C$ with $R(K) = C(K),$ but $P^t(1,K)$ still has bounded point evaluations. The first part of this paper is to extend the above result of \cite{bm} to $P^t(1, \phi_1,...,\phi_N,K)$ and $HP(\phi_1,...,\phi_N,K)$ for some smooth functions $\phi_1,...,\phi_N.$ The following theorem, which connects analytic capacity with special types of bounded point evaluation estimations, is essential for our main results. Notice that Theorem \ref{MTheorem1} (1) extends Lemma B in \cite{ars}.

\begin{Theorem} \label{MTheorem1}
There exist absolute constants $\epsilon_N, C_N > 0$ that only depend on $N$. If 
\[
 \ \gamma (B (\lambda_0, \delta) \setminus K) < \epsilon_N \delta, \tag{1-3}
 \]
then

(1) 
\[
 \ |p_N(\lambda)| \le \dfrac{C_N}{\delta ^{N+2}}\int _{K\cap\bar B (\lambda_0, \delta)} \left |\sum_{k=0}^N p_k \bar z ^k \right | dA,~\lambda \in  \bar B \left (\lambda_0, \frac{1}{2}\delta \right ), \tag{1-4}
 \]
 
(2)
\[
 \ |p_N(\lambda)| \le \dfrac{C_N}{\delta ^N} \left \|\sum_{k=1}^N p_k \bar z ^k + r \right \|_{K\cap\bar B (\lambda_0, \delta)}, ~\lambda \in  \bar B \left (\lambda_0, \frac{1}{2}\delta \right ), \tag{1-5}
 \]
where $r\in Rat(K\cap\bar B (\lambda_0, \delta))$ and $p_k\in\mathcal P.$
 \end{Theorem}

The proof of Theorem \ref{MTheorem1} depends on a careful modification of Thomson's coloring scheme on dyadic squares.  Thomson's coloring scheme, for a point $a \in \mathbb C$ and a
positive integer $m,$ starts with a dyadic square of side length $2^{-m}$ containing $a$
and either terminates at some finite stage or produces an infinite sequence of
annuli surrounding $a.$ These annuli are made up of dyadic squares colored red (heavy square). When the scheme terminates at some finite stage, one can find a path consisting of many dyadic squares colored green (light square). The definition of a light square in \cite{thomson}
(see also page 168 of \cite{thomson3} or page 461 of \cite{ars}) only works for $P^t(1,K)$ and $HP(\bar z, K).$ Our definition of a light square (2-2) allows us to recursively extend (1-4) and (1-5) for $N > 1.$

Define
 \[
 \ BA(\lambda_0, \delta ) = \{f\in C(\mathbb C_\infty): ~\exists n,~f\in C^{(n)}(\bar B (\lambda_0, \delta )),~ \bar\partial^nf |_{\bar B (\lambda_0, \delta )} = 0\},
 \]
where $\bar\partial$ is the Cauchy-Riemann operator.
We now state our first main result.

\begin{Theorem} \label{BPETheorem}
Let $\lambda_0 \in K$ be a nonpeak point for $R(K)$ and $1\le t < \infty .$ If there exist $\delta > 0$ and $F_1, F_2,...,F_N\in BA(\lambda_0, \delta )$ such that
 \[
 \ det[ \bar\partial^i F_j (\lambda _0) ]_{n\times n} \ne 0,
 \] 
then

(1) $\lambda _0$ is an analytic bounded point evaluation for $P^t(1, F_1,...,F_N, K).$ In particular, $\lambda _0$ is an analytic bounded point evaluation for $P^t(1, \bar z,...,\bar z^N, K).$

(2) $\lambda _0$ is an analytic bounded point evaluation for $HP(F_1,...,F_N, K).$ In particular, $\lambda _0$ is an analytic bounded point evaluation for $HP(\bar z,...,\bar z^N, K).$
\end{Theorem}

Let $\Lambda$ be a constant coefficient elliptic differential operator in $R^2$. For a compact $K \subset \mathbb C,$ let $H (K, \Lambda)$ and $h (K, \Lambda)$ denote the uniform closures in $C(K)$ of the set
 \[
 \ \{f\in K : \Lambda f = 0 ~ \text{in some neighborhood of } K\}
 \]
and the set
 \[
 \ C(K) \cap\{f\in K : \Lambda f = 0 ~ \text{in the interior of } K\}
 \]
respectively. Notice that, for $\Lambda  =  \bar\partial,$ the space $H (K, \bar\partial) = R(K)$ and $h (K, \bar\partial) = A(K).$ For $\Lambda  =  \bar\partial ^ n,$ the $nth$ power of Cauchy-Riemann operator, the space 
 \[ 
 \ H (K, \bar\partial^n) = clos ( R(K) + \bar z R(K) + ...  + \bar z ^{n-1 } R(K))
 \]
and 
 \[
 \ h (K, \bar\partial^n) = clos( A(K) + \bar z A(K) + ...  + \bar z ^{n-1 } A(K)).
 \] 
One of uniform approximation problems is the following:

\begin{Problem}\label{PP}
Find necessary and sufficient conditions for $K$ so that 
$H (K, \Lambda) = h (K, \Lambda).$
\end{Problem}

A complete solution for $\Lambda = \Delta$ was obtained by \cite{Deny} and \cite{Keldysh} using a duality argument relying on classical potential theory. Let $Cap$ denote the Wiener capacity in potential theory. Deny and Keldysh show that the identity $H (K, \Delta) = h (K, \Delta)$ occurs if and only if for each open ball $B$ one has $Cap(B \setminus int K) = Cap(B \setminus K).$
 
Using a constructive scheme for uniform approximation (based on a localization operator), \cite{Vitushkin} proves that the identity $H (K, \bar\partial) = h (K, \bar\partial)$ occurs if and only if for each open disc $O$ one has $\alpha(O \setminus int K) = \alpha(O \setminus K).$

The inner boundary of $K$, denoted by $\partial_i K$, is the set of boundary points which do not belong to the boundary of any connected component of $\mathbb C\setminus K$. The remarkable paper \cite{Tol04} proves that the continuous analytic capacity is semiadditive. The result implies an affirmative answer to the so called “inner boundary conjecture” (see \cite{vm}, Conjecture 2). That is, if $\alpha(\partial_i K) = 0,$ then $R(K) = A(K).$ 

For $\Lambda = \bar\partial^2,$ \cite{TrentWang} show if K is a compact subset without interior, then $H (K, \bar\partial^2) = h (K, \bar\partial^2) = C(K).$ \cite{ver} proves that each Dinicontinuous function in $h (K, \bar\partial^2)$ belongs to $H (K, \bar\partial^2).$ Finally, the excellent paper \cite{mazalov} completely solved the problem by proving $H (K, \bar\partial^2) = h (K, \bar\partial^2)$ for any compact subset $K.$

In \cite{bcf}, the authors consider an interesting analogous problem: find necessary and sufficient conditions so that $P(K) + P(K) \bar z^n$ is dense in $A(K) + A(K) \bar z^n.$ The paper obtained some results for a Carath´eodory compact set $K$ with $n \ge 2$ (see \cite{bcf}, Theorem 1).

We define
 \[
 \ hP(\phi_1,\phi_2,...,\phi_N, K) = clos \left ( \sum_{i=1}^N P(K) \phi_i + A(K) \right).
 \]
As analogous to Problem \ref{PP}, we are interested in the following problem: 

\begin{Problem}\label{MP}
Find necessary and sufficient conditions so that 
 \[
 \ HP(\phi_1,\phi_2,...,\phi_N, K) = hP(\phi_1,\phi_2,...,\phi_N, K).\tag{1-6}
 \]
\end{Problem}

For $N = 1,$ the problem was studied by several authors. In \cite{thomson3}, Thomson proves if $R(K) \ne C(K),$ then $HP(\bar z, K)$ is not equal to $C(K).$ \cite{yang1} and \cite{yang2} study the generalized space $HP(g, K)$ and prove that for a smooth function $g$ with $\bar \partial g \ne 0,$ then $HP(g, K) = hP(g, K)$ if and only if $ A(K) = R(K).$

In the second part of this paper, we will study Problem \ref{MP} when $N > 1.$ Our second main theorem is the following: 

\begin{Theorem}\label{MTheorem2}
Let $q_1(z,\bar z),q_2(z,\bar z),...,q_N(z,\bar z))$ be polynomials in two variables $z$ and $\bar z.$ 
If
 \[
 \ A(K) \subset HP(q_1,...,q_N,K),
 \]
then $A(K) = R(K).$ In particular, if $R(K) \ne C(K),$ then $HP(q_1,...,q_N,K) \ne C(K).$ Notice that an important special case is that $q_1(z,\bar z) = \bar z,$ $q_2(z,\bar z) = \bar z^2,$ ...,$q_N(z,\bar z) = \bar z^N.$
\end{Theorem}

By Stone-Weierstrass theorem, $ \sum_{k=1}^{\infty} P(K) \bar z ^k + R(K) $ is dense in $C(K).$ So it is critical here to assume that $N$ is a finite integer.

In \cite{bcf}, the authors are interested in the question of finding necessary and sufficient conditions so that 
 \[
 \ clos \left (A(K) + \sum_{i=1}^N A(K)\bar z^{d_i}\right ) = clos \left (P(K) + \sum_{i=1}^N P(K)\bar z^{d_i}\right ), \tag{1-7}
 \] 
where $d_1,...,d_N$ are positive integers. Our theorem implies if (1-7) holds, then $A(K) = R(K).$ So (1-7) is equivalent to 
 \[
 \ clos \left (R(K) + \sum_{i=1}^N R(K)\bar z^{d_i}\right ) = clos \left (P(K) + \sum_{i=1}^N P(K)\bar z^{d_i}\right ).
 \]

The following result shows that Theorem \ref{MTheorem2} will not hold if we replace $q_n$ by certain functions.

\begin{Proposition}\label{Prop}
There exist a compact subset $K\subset \mathbb C$ and a function $\phi\in A(K)$ such that $R(K) \ne A(K),$ but $A(K) = HP(\phi , K).$
\end{Proposition}  

The proposition raises the following question:

\begin{Problem} For a compact subset $K$ of $\mathbb C,$ is there a function $\phi\in A(K)$ such that $A(K) = HP(\phi, K)?$
\end{Problem}

It seems that one might be able to use \cite{Tol04} characterization of continuous analytic capacity to find a proper finite Borel measure $\nu$ supported on the inner boundary of $K$ such that $A(K) = HP(\mathcal C\nu, K).$

We prove Theorem \ref{MTheorem1} and Theorem \ref{BPETheorem} in Section 2. In section 3, we prove Theorem \ref{MTheorem2} and construct a compact subset $K$ to prove Proposition \ref{Prop}.

\section{Bounded Point Evaluations} 

In this section, we will prove Theorem \ref{MTheorem1} and Theorem \ref{BPETheorem}.
\cite{Tol03} proves the astounding result about analytic capacity $\gamma$ that 
implies the semiadditivity of analytic capacity. That is,
\[
\ \gamma (\bigcup_{i = 1}^m E_i) \le C_T \sum_{i=1}^m \gamma(E_i) \tag{2-1}
\]
where $C_T$ is an absolute constant.

For a square $S$ (also denoted by $S(c_S, d_S)$), whose edges are parallel to x-axis and y-axis, let $c_S$ denote the center and $d_S$ denote the side length. For $a>0,$ $aS$ is a square with the same center of $S$ ($c_{aS} = c_S$) and the side length $d_{aS} = ad_S.$ For a given $\epsilon > 0,$ a closed square $S$ is defined to be light $\epsilon$ if
 \[
 \ \gamma (Int(S) \setminus K)) > \epsilon  d_S.\tag{2-2}
 \]
A square is called heavy $\epsilon$ if it is not light $\epsilon$.
Let $R = \{ z: -1/2 < Re(z),Im(z) < 1/2 \}$ and $Q = \bar{\mathbb D}\setminus R.$ 

We now sketch our version of
Thomson's coloring scheme for $Q$ with a given $\epsilon$ and a positive integer $m.$ We refer the reader to \cite{thomson} and \cite{thomson3} section 2 for details.
 
For each integer $k > 3$ let $\{S_{kj}\}$ be an enumeration of the closed squares contained in $\mathbb C$ with edges of length $2^{-k}$
parallel to the coordinate axes, and corners at the points whose coordinates
are both integral multiples of $2^{-k}$ (except the starting square $S_{m1}$, see (3) below). 
In fact, Thomson'scoloring scheme is just needed to be modified slightly as the following:

(1) Use our definition of a light $\epsilon$ square (2-2).

(2) A path to $\infty$ means a path to any point that is outside of $Q$ (replacing the polynomially convex hull of $\Phi$ by $Q$).

(3) The starting yellow square $S_{m1}$ in the $m$-th generation is $R.$ Notice that the length of $S_{m1}$ in $m$-th generation is $1$ (not $2^{-m}$).

We will borrow notations that are used in Thomson's coloring scheme such as $\{\gamma_n\}_{n\ge m}$ and $\{\Gamma_n\}_{n\ge m},$ etc. We denote 
 \[
 \ YellowBuffer_m = \sum _{k = m+1}^\infty k^2 2^{-k}.
 \]
Two things can happen (depending on $m$):

Case I. The scheme terminates, in our setup, this means Thomson's coloring scheme reaches a square $S$ in $n$-th generation that is not contained in $Q.$ One can construct a polygonal path $P,$ which connects the centers of adjacent squares, from the center of a square (contained in $Q$) adjacent to $S$ to the center of a square adjacent to $R$ so that the orange (non green in Thomson's coloring scheme) part of length is no more than $YellowBuffer_m.$ Let $GP = \cup S_j,$ where $\{S_j\}$ are all light $\epsilon$ squares with  $P\cap S_j \ne \emptyset .$ By Tolsa's Theorem (2-1), we see
 \[
 \ \gamma (P) \le C_T (\gamma (Int(GP)) + YellowBuffer_m).
 \]
Since $P$ is a connected set, $\gamma (P) \ge \frac{0.1}{4}$ (Theorem 2.1 on page 199 of \cite{gamelin}).  We can choose $m $ to be large enough so that
 \[
 \ \gamma (Int(GP)) \ge \dfrac{1}{40C_T} - YellowBuffer_{m} = \epsilon_m > 0. \tag{2-3}
 \]

Case II. The scheme does not terminate. In this case, one can construct a sequence of heavy $\epsilon$ barriers inside $Q,$ that is, $\{\gamma_n\}_{n\ge m}$ and $\{\Gamma_n\}_{n\ge m}$ are infinite.

For simplicity, we will use $scheme(Q,\epsilon,m, \gamma_n, \Gamma_n, n\ge m)$ to stand for our version of Thomson's coloring scheme.

If a function $f$ is analytic at $\infty,$ then $f$ can be represented by its Laurent series
 \[
 \ f(z) = f(\infty ) + a_1 (z- z_0)^{-1} + a_2(z - z_0)^{-2} + ...
 \]
in a neighborhood of infinite. We define $f'(\infty )$ to be $a_1$ and $\beta (f, z_0)$ to be
$a_2.$ The number $f'(\infty )$ does not depend on $z_0,$ but $\beta (f, z_0)$ does depend
on the choice of $z_0.$

\begin{Lemma}\label{prop1} For a square $T,$ if 
 \[
 \ \gamma(Int(T) \setminus K) \ge \epsilon_1 d_T,
\]
then for two complex numbers $|\alpha| \le 1$ and $|\beta| \le 1,$ there exists a function $f$ in $C(\mathbb C_\infty)$ such that the following hold:

(1) $\|f\| \le \frac{54}{\epsilon_1^3};$

(2) $f\in R(\mathbb C_\infty \setminus (Int(T)\setminus K))$ ;

(3) $f(\infty) = 0;$

(4) $f'(\infty ) = \alpha d_T;$

(5) $\beta (f,c_T) = \beta d_T^2.$
\end{Lemma}

Proof: There exists a function $f_1$ in $C(\mathbb C_\infty)$ such that $\|f_1\| \le 1,$
$f_1$ is analytic off a compact subset of $Int(T) \setminus K,$ $f_1(\infty) = 0,$ and $f_1'(\infty ) > \epsilon _1 d_T/2.$ Set $f_2 = d_Tf_1/f_1'(\infty).$ Then by Theorem 2.5 of \cite{gamelin} on page 201, we get
 \[
 \ |\beta (f_2,0)|  \le \dfrac{12}{\epsilon _1} d_T^2.
 \] 
Let
 \[
 \ f = \alpha \left (f_2 - \dfrac{\beta (f_2,0)}{d_T^2} f_2^2\right ) + \beta f_2^2,
 \]
then $f\in R(\mathbb C_\infty \setminus (Int(T)\setminus K))$ and satisfies the conditions (1)-(5).

Let $\varphi$ be a smooth function with compact support. The localization operator
$T_\varphi$ is defined by
 \[
 \ (T_\varphi f)(\lambda) = \dfrac{1}{\pi}\int \dfrac{f(z) - f(\lambda)}{z - \lambda} \bar\partial \varphi (z) dA(z),
 \]
where $f$ is a continuous function on $\mathbb C_\infty.$ One can easily prove the following
norm estimation for $T_\varphi:$
 \[
 \ \| T_\varphi f\| \le  4\|f\| diameter(supp \varphi)\|\bar\partial \varphi\|.
 \]

\begin{Lemma} \label{LightRoute}
Suppose Case I of $scheme(Q,\epsilon,m, \gamma_n, \Gamma_n, n\ge m)$ is true, then
 \[
 \ \gamma (\mathbb D \setminus K) \ge \epsilon_1, \tag{2-4}
 \]
where
 \[
 \ \epsilon_1 = 10^{-8}\epsilon ^3 \epsilon _m \tag{2-5}
 \]
and $\epsilon _m$ is from (2-3).
\end{Lemma}

Proof: we will follow the second part of proof of Lemma B in \cite{ars} on pages 464-465 with slight modifications. Let $GP = \cup S_j,$ where $\{S_j\}$ are light $\epsilon$ squares discussed above, so that 
$\gamma (Int(GP)) \ge \epsilon_m.$ For each j let
$z_j$ be the center of $S_j ,$
$d_j$ be the edge length of $S_j ,$
$Q_j,$ $R_j$ be the closed squares with center $z_j$ and sides parallel to the
coordinate axes of lengths $\frac{7}{6} d_j = \delta _j ,$ $\frac{2}{3}d_j$ respectively.
The collection $\{S_j\}$ has the following properties (see (2.16)-(2.18) on page 464 of \cite{ars}):

(a) No point lies in more than four $Q_j$'s.

(b) There are $C^\infty$ functions $\phi_j$ with $0 \le \phi_j \le 1,$ $spt (\phi_j) \subset Int(Q_j) ,$
$\|\bar \partial \phi_j \| \le \frac{50}{d_j} ,$ $\phi_j = 1$ on $R_j ,$ and $ \sum \phi_j = 1$ on $GP.$

(c) For each $z\in \mathbb C,$
 \[
 \ \sum \min \left \{1, \dfrac{\delta _j^3}{|z - z_j|^3} \right \} \le 10,000.
 \]
Let $f\in C(\mathbb C^\infty )$ such that $f$ is analytic off a compact subset of $Int(GP),$ $f(\infty ) = 0,$ $\|f\|_\infty = 1,$ and $f'(\infty ) > \frac{\epsilon_m}{2}.$ Then from (b), we see that $f - \sum T_{\phi_j}f$ is zero on $GP$ and analytic off $GP.$ Hence,
 \[
 \  f(z) = \sum T_{\phi_j}f(z)
 \]
for all $z\in \mathbb C$ and
 \[
 \ T_{\phi_j}f(\infty ) = 0, ~\|T_{\phi_j}f\| \le 400, ~ |(T_{\phi_j}f)'(\infty )| \le 400d_j, ~ \text{and}~ |\beta (T_{\phi_j}f, z_j)| \le 400 d_j^2. 
 \]
For $\alpha = \frac{(T_{\phi_j}f)'(\infty )}{400d_j}$ and $\beta = \frac{\beta (T_{\phi_j}f, z_j)}{400 d_j^2},$ using Lemma \ref{prop1} for the light $\epsilon$ square $S_j,$ we find a function $f_j$ in $C(\mathbb C_\infty)$ such that $\|f_j\| \le \frac{54}{\epsilon^3},$  $f_j \in R(\mathbb C_\infty \setminus (Int(S_j)\setminus K)),$ and  $\frac{T_{\phi_j}f}{400} - f_j$ has triple zeros at $\infty.$ Therefore, $(\frac{T_{\phi_j}f}{400} - f_j)(z-z_j)^3$ is analytic on $\mathbb C_\infty \setminus Q_j.$ Using the Maximum Modulus
Theorem, we see if $z \in \mathbb C_\infty \setminus Q_j ,$ then 
 \[
 \ \left | \frac{T_{\phi_j}f(z)}{400} - f_j(z) \right | \le (1 + \dfrac{54}{\epsilon ^3})\dfrac{\delta_j^3}{|z-z_j|^3} \le \dfrac{55}{\epsilon ^3}\dfrac{\delta_j^3}{|z-z_j|^3}.
 \]
Hence, for 
$z \in \mathbb C_\infty ,$  
 \[
 \ \left | \frac{T_{\phi_j}f(z)}{400} - f_j(z) \right | \le \dfrac{55}{\epsilon ^3} \min \left (1, \dfrac{\delta_j^3}{|z-z_j|^3} \right ).
 \]
Set $F = 400\sum f_j,$ then $F$ is analytic off a compact subset of $Int(GP)\setminus K,$ $F(\infty ) = 0,$ $F'(\infty ) = f'(\infty ),$ and
 \[
 \ \begin{aligned}
 \ \|F \| _\infty \le & \|f \| _\infty + \sum \|T_{\phi_j}f - 400f_j\|_\infty \\
 \ \le & 1 + \dfrac{2200}{\epsilon ^3} \sum \min \left (1, \dfrac{\delta_j^3}{|z-z_j|^3} \right ) \\ 
 \ \le & \dfrac{5\times10^7}{\epsilon ^3}.
 \ \end{aligned}
 \]
where the last step follows from (c). Therefore,
 \[
 \ \gamma (\mathbb D \setminus K) \ge \gamma (Int(GP) \setminus K) \ge \dfrac{F'(\infty )}{\frac{5\times10^7}{\epsilon ^3}} > 10^{-8}\epsilon ^3 \epsilon _m = \epsilon _1.
 \]

\begin{Lemma} \label{HBarrier}
Suppose Case II of $scheme(Q,\epsilon,m, \gamma_n, \Gamma_n, n\ge m)$ is true. 

(1) If there exists $\epsilon _N > 0$ such that for every heavy $\epsilon$ square $S$ in $scheme(Q,\epsilon,m, \gamma_n, \Gamma_n, n\ge m),$ 
 \[
 \ dist \left (\bar z ^N, P^1 (1, \bar z, ..., \bar z^{N-1}, K \cap S ) \right ) \ge \epsilon_N d_S^{N+2}, \tag{2-6}
 \] 
then there exists a constant $C_{m,N}$ (depends on $m$ and $N$) such that
 \[
 \ |p_N(\lambda )| \le C_{m,N}\left \| \sum_{j=0}^N p_j \bar z^j \right\|_{L^1(K\cap Q)}, ~ \lambda \in R, \tag{2-7}
 \]
where $p_j \in \mathcal P$ and $N = 1,2,....$ For $N = 0,$ if (2-6) is replaced by 
 \[
 \ Area(K \cap S) \ge \epsilon_0 d_S^{2}, \tag{2-8}
 \]
then (2-7) holds for $N=0.$

(2) If there exists $\epsilon _N > 0$ such that for every heavy $\epsilon$ square $S$ in $scheme(Q,\epsilon,m, \gamma_n, \Gamma_n, n\ge m),$ 
 \[
 \ dist \left(\bar z^N, HP(\bar z,\bar z^2,...,\bar z^{N-1}, K\cap S) \right ) \ge \epsilon_N d_S ^N, \tag{2-9}
 \] 
then there exists a constant $C_{m,N}$ (depends on $m$ and $N$) such that
 \[
 \ |p_N(\lambda )| \le C_{m,N}\left \| \sum_{j=1}^N p_j \bar z^j + r \right\|_{K\cap Q}, ~ \lambda \in R, \tag{2-10}
 \]
where $r \in Rat(K\cap Q),$ $p_j\in \mathcal P,$   and $N = 1,2,....$
\end{Lemma}

Proof: The proofs of section 4 in \cite{thomson} and \cite{thomson3} will work if we make the following modifications:

(a) For $w\in \gamma_n \cap S,$ where $S$ is a heavy $\epsilon$ square, by the Hahn-Banach theorem,  there exists a finite Borel measure $\sigma_w$ supported in $K\cap S$ and $\|\sigma_w\| = 1$ such that for (1), $\sigma_w = \Phi_w dA$ with $\Phi_w \in L^\infty (K\cap S),$
(2-6) becomes
\[
 \ dist \left (\bar z ^N, P^1 (1, \bar z, ..., \bar z^{N-1}, K \cap S ) \right ) = \int \bar z ^N d \sigma _w \ge \epsilon_N d_S^{N+2}, 
 \]
where $\int f d \sigma _w = 0$ for $f \in P^1 (1, \bar z, ..., \bar z^{N-1}, K \cap S ),$  and (2-8) becomes
 \[
 \  \sigma_w = \chi_{K\cap S} dA,~ Area(K \cap S) = \int d \sigma _w \ge \epsilon_0 d_S^{2};
 \]
 and for (2), (2-9) becomes
 \[
 \ dist \left(\bar z^N, HP(\bar z,\bar z^2,...,\bar z^{N-1}, K\cap S) \right ) = \int \bar z ^N d \sigma _w  \ge \epsilon_N d_S ^N,
 \]
where $\int f d \sigma _w = 0$ for $f \in HP(\bar z,\bar z^2,...,\bar z^{N-1}, K\cap S). $
Define $\tau_w = \frac{\bar z ^N d \sigma _w}{\int \bar z ^N d \sigma _w}.$

(b) Set $L = e_\lambda$ ($e_\lambda (f) = f(\lambda )$) for $\lambda \in R.$ Use the same argument as in section 4 in \cite{thomson3}, one can construct $\mu_{n+q}$ defined by linear combination of $\tau_w,$ then
 \[
 \ \| \mu _{n+q} \| \le \epsilon_N^{-q-1}2^{(N+2)(n+q)} 4^q ((n+q-1)...n)^{-2}
 \]
for (1) and
 \[
 \ \| \mu _{n+q} \| \le \epsilon_N^{-q-1}2^{N(n+q)} 4^q ((n+q-1)...n)^{-2}
 \]
for (2).
Let $\mu = \mu_n +...+ \mu_{n+q}+ ...,$ then $\| \mu \| \le C_{m,N},$ 
for (1)
 \[
 \ p_N(\lambda ) = \int \sum_{j=0}^N p_j \bar z^j d \mu,
 \]
and for (2)
 \[
 \ p_N(\lambda ) = \int (\sum_{j=1}^N p_j \bar z^j +r ) d \mu.
 \]
Clearly, the support of $\mu$ is outside $R.$ The proof is completed.

The idea to prove our Theorem \ref{MTheorem1} is to find sufficient small $\epsilon$ so that Case II of $scheme(Q,\epsilon,m, \gamma_n, \Gamma_n, n\ge m)$ is true. Then we use mathematical induction to show that for every heavy $\epsilon$ square, we can find $\epsilon_N$
so that (2-6), (2-8), and (2-9) all hold. We will demonstrate the idea by proving the following corollary, which is the case for $N=0$ in (1) of Theorem \ref{MTheorem1} and Lemma B in \cite{ars}. 

\begin{Corollary} \label{corollaryARS}
There are absolute constants $\epsilon _0 > 0$ and $C_0 < \infty$ with the
following property. For a compact subset $K\subset \mathbb C,$ let $R > 0$ and $\gamma (R\mathbb D \setminus K) < \epsilon_0 R.$ Then 
\[
\ |p(\lambda )| \le \dfrac{C_0}{R^2} \int _{( R\mathbb D)\cap K} |p| \frac{dA}{\pi}
\]
for $|\lambda | \le \frac{R}{2}$ and all $p \in \mathcal P.$
\end{Corollary}

Proof: Since $\gamma (R\mathbb D \setminus K) = R\gamma (\mathbb D \setminus \frac{K}{R}),$ by a simple changing of variables from $z$ to $Rz,$ we assume $R=1.$
Let $\epsilon _1 = \epsilon_0 = 10^{-8}\epsilon ^3 \epsilon _m$ in (2-5) with $\epsilon < \sqrt{\frac{1}{4\pi}}.$ Then from Lemma \ref{LightRoute}, we conclude that Case II of $scheme(Q,\epsilon,m, \gamma_n, \Gamma_n, n\ge m)$ must be true. Let $S$ be a heavy $\epsilon$ square, then $\gamma (Int(S)\setminus K) \le \epsilon d_S.$ By Theorem 3.2 on page 204 of \cite{gamelin}, we get
 \[
 \ Area(S\setminus K ) \le 4\pi \gamma (Int(S)\setminus K)  ^2 \le 4\pi \epsilon ^2 d_S^2. 
 \]  
Therefore,
 \[
 \ Area (S \cap K) \ge (1 - 4\pi \epsilon ^2) d_S^2.
 \]
So (2-8) holds and the corollary follows from Lemma \ref{HBarrier}.

Let $\phi$ be a smooth function supported in $\mathbb D$  such that:
 \[
 \ 0\le \phi \le 1,~ \phi (z) = \phi (|z|),~ \|\bar\partial ^{N} \phi \|<C_N^\phi,~ \int \phi dA = 1.\tag{2-11}
 \]

Proof of Theorem \ref{MTheorem1} (1): We only need to prove the case that $\lambda _0 = 0$ and $\delta = 1.$ In fact, using the elementary properties of analytic capacity (see p. 196 of \cite{gamelin}), one sees that condition (1-3) is equivalent to
 \[
 \ \gamma \left (B (0, 1) \setminus \dfrac{K -\lambda_0}{\delta} \right ) < \epsilon_N.
 \]
The inequality (1-4)
 \[
 \begin{aligned}
 \ |p_N(\lambda)| & \le \dfrac{C_N}{\delta ^{N+2}} \left \|\sum_{k=0}^N p_k \bar z ^k \right \|_{L^1(K\cap\bar B (\lambda_0, \delta))} \\
 \ & =   \dfrac{C_N}{\delta ^2} \left \|\sum_{k=0}^N q_k (\dfrac{\bar z - \bar \lambda _0}{\delta}) ^k\right \|_{L^1(K\cap\bar B (\lambda_0, \delta))}, 
 \end{aligned}
 \]
for $\lambda \in  \bar B (\lambda_0, \frac{1}{2}\delta ),$ where $q_N = p_N,$ $q_k (1 \le k \le N-1)$ are certain linear combinations of $p_k,$ is equivalent to 
 \[
 \ |p_{0N}(\lambda)| \le  C_N \left \|\sum_{k=0}^N p_{0k} \bar z  ^k \right \|_{L^1(\frac{K -\lambda_0}{\delta}\cap\bar B (0, 1))}, ~\lambda \in  \bar B \left (0, \frac{1}{2} \right ),
\]
where $p_{0k}(z) = q_k(\delta z + \lambda _0).$
We will assume that $\lambda _0 =0$ and $\delta = 1$ in the rest of the proof.

We use mathematical induction for $N.$ The case $N=0$ is directly implied by Corollary \ref{corollaryARS}.

Now we assume that Theorem \ref{MTheorem1} (1) holds for $k =0,1,2,...,N.$ Set
 \[
 \ \epsilon = min (c_1, \dfrac{\epsilon _{N}}{2}, \dfrac{\epsilon _{N-1}}{2^2},...,\dfrac{\epsilon _1}{2^{N}},\dfrac{\epsilon _0}{2^{N+1}}),
 \] 
where $c_1 > 0$ will be determined later, and $\epsilon _{N+1}$ as in (2-5), that is, $\epsilon _{N+1} = 10^{-8}\epsilon ^3 \epsilon _m.$ Since $\gamma (B (0, 1) \setminus K) < \epsilon_{N+1}$ (assumption (1-3)), from Lemma \ref{LightRoute}, we conclude that Case II of $scheme(Q,\epsilon,m, \gamma_n, \Gamma_n, n\ge m)$ is true. For every heavy $\epsilon$ square $S,$ we need to prove (2-6) holds. 

We assume
 \[
 \ dist \left (\bar z ^{N+1}, P^1 (1, \bar z, ..., \bar z^N, K \cap S ) \right ) < \frac{1}{2} d_S^{N+3}, 
 \]
otherwise (2-6) already holds. Without loss of generality, we assume the center of $S$ is zero. There are polynomials $p_0,p_1,..., p_N$ such that
 \[
 \ \left \| \bar z^{N+1} + \sum _{j=0}^N p_j \bar z^j \right \| _{L^1(K\cap S)} \le 2 dist \left (\bar z ^{N+1}, P^1 (1, \bar z, ..., \bar z^N, K \cap S ) \right ) <d_S^{N+3}.\tag{2-12}
 \]
Hence,
 \[
 \ \left \|\sum _{j=0}^N p_j \bar z^j \right \| _{L^1(K\cap S)} \le d_S^{N+3} +  \|\bar z^{N+1}\|_{L^1(K\cap S)} \le 2d_S^{N+3}.
 \]
Since $S$ is a heavy $\epsilon$ square, we get
 \[
 \ \gamma (B(0, \frac{d_S}{2} ) \setminus K) \le \gamma (S \setminus K) < \epsilon d_S \le \epsilon_N \frac{d_S}{2}.
 \]
Applying the induction assumption for $N,$ from (1-4), we have
 \[
 \ |p_N(\lambda)| \le \dfrac{C_N}{d_S^{N+2}} \left \|\sum _{j=0}^N p_j \bar z^j  \right \| _{L^1(K\cap \bar B(0, \frac{d_S}{2}))} < 2C_N d_S,
 \]
for $\lambda \in \bar B(0, \frac{d_S}{4}).$ This implies
 \[
 \ \|\bar z^{N+1} + p_N(z)\bar z^N\|_{L^1(K\cap \bar B(0, \frac{d_S}{4}))}\le 3C_Nd_S^{N+3} .
 \]
From (2-12), we conclude
 \[
 \ \left \|\sum _{j=0}^{N-1} p_j \bar z^j \right \| _{L^1(K\cap \bar B(0, \frac{d_S}{4}))} \le 4C_Nd_S^{N+3}.
 \]
In general, there is an absolute constant $C_{N+1} >0 $ so that
 \[
 \ \left \|\sum _{j=0}^{N-k} p_j \bar z^j\right \| _{L^1(K\cap \bar B(0, \frac{d_S}{2^{k+1}}))} \le C_{N+1}d_S^{N+3}.
 \]
Since $S$ is a heavy $\epsilon$ square, we get
 \[
 \ \gamma (B(0, \frac{d_S}{2^{k+1}} ) \setminus K) \le \gamma (S \setminus K) < \epsilon d_S \le \epsilon_{N-k} \frac{d_S}{2^{k+1}}.
 \]
Apply the induction assumption for $N-k$, from (1-4), we get
 \[
 \ |p_{N-k}(\lambda)| \le \dfrac{C_{N-k}}{d_S^{N-k+2} }\left \|\sum _{j=0}^{N-k} p_j \bar z^j \right \| _{L^1(K\cap \bar B(0, \frac{d_S}{2^{k+1}}))} \le C_{N-k} C_{N+1} d_S^{k+1},
 \]
for $\lambda \in \bar B(0, \frac{d_S}{2^{k+2}}).$
So there is an absolute constant which we still use $C_{N+1}> 0$ such that 
 \[
 \ \left | \bar z^{N+1} + \sum _{j=0}^N p_j (z) \bar z^j \right | \le C_{N+1}d_S^{N+1}
 \]
for $z\in \bar B(0, \frac{d_S}{2^{N+1}}).$
Let $g_{N+1} = \bar z^{N+1} + \sum _{j=0}^N p_j \bar z^j,$ then $|g_{N+1} (z)| \le C_{N+1}d_S^{N+1}$ on $z\in \bar B(0, \frac{d_S}{2^{N+1}}).$

Let $\phi_{N+1}(z) = \phi (\frac{2^{N+1}z}{d_S}),$ where $\phi$ is in (2-11), then $spt(\phi_{N+1}) \subset B(0, \frac{d_S}{2^{N+1}}),$
 \[
 \ 0\le \phi_{N+1} \le 1,~ \|\bar\partial ^{N+1} \phi_{N+1} \| < C_{N+1}^{\phi_{N+1}}/d_S^{N+1},~ \int \phi_{N+1} dA = \frac{d_S^2}{4^{N+1}}.
 \]

Then
 \[
 \ \begin{aligned}
 \ & (N+1)! \frac{d_S^2}{4^{N+1}} \\
 \ = & (N+1)! \int \phi_{N+1} dA \\
 \ = & \left | \int g_{N+1} \bar\partial ^{N+1}  \phi_{N+1} dA  \right | \\ 
 \ \le & \dfrac{C_{N+1}^{\phi_{N+1}}}{d_S^{N+1}} \left (\left \| g_{N+1}  \right \| _{L^1(K\cap \bar B(0, \frac{d_S}{2^{N+1}}))} + \left \| g_{N+1}  \right \| _{L^1(\bar B(0, \frac{d_S}{2^{N+1}}) \setminus K)} \right )\\
 \ \le & \dfrac{C_{N+1}^{\phi_{N+1}}}{d_S^{N+1}}\left \| g_{N+1} \right \| _{L^1(K\cap S)}+ C_{N+1}C_{N+1}^{\phi_{N+1}} Area( B(0, \frac{d_S}{2^{N+1}}) \setminus K) \\
 \ \le & 2 \dfrac{C_{N+1}^{\phi_{N+1}}}{d_S^{N+1}}dist \left (\bar z ^{N+1}, P^1 (1, \bar z, ..., \bar z^N, K \cap S ) \right )+ 4\pi C_{N+1}C_{N+1}^\phi \gamma ( Int(S) \setminus K)^2
 \ \end{aligned}
 \]
where the last step follows from (2-12) and Theorem 3.2 on page 204 of \cite{gamelin}. Now choose
 \[
 \ c_1^2 = \dfrac{(N+1)!}{2^{2N+5}\pi C_{N+1}C_{N+1}^{\phi_{N+1}}},
 \]
then since $S$ is a heavy $\epsilon$ square, we have
 \[
 \ dist \left (\bar z ^{N+1}, P^1 (1, \bar z, ..., \bar z^N, K \cap S ) \right ) \ge \dfrac{(N+1)!}{2^{2N+4}C_{N+1}^{\phi_{N+1}}}d_S^{N+3}.
 \]
So (2-6) holds and the theorem follows from Lemma \ref{HBarrier} (1).

For a smooth function $\varphi$ with compact support, $T_\varphi$ is a bounded linear operator on $C(K).$ Let $M$ be the space of finite complex Borel measures supported on $K$ ( = $C(K)^*$), then $T_\varphi^*$ is a bounded linear operator on $M.$ Moreover, for $\mu\in M,$
 \[
 \ \int T_\varphi f d\mu = \int f dT_\varphi^*\mu
 \]
and
 \[
 \ \|T_\varphi^*\mu \| \le 4 diameter(supp\varphi)\|\bar\partial\varphi\| \|\mu\|.   \tag{2-13}
 \]
Consequently, $T_\varphi^*\mu \perp R(K)$ for each $\mu \perp R(K).$  

\begin{Lemma}\label{propN} Let $S$ be a square with center $0$ and $d_S < 1.$ Let $g_N = \bar z^N + \sum_{k=1}^{N-1} p_k \bar z^k,$ where $p_k$ are polynomials and $\|g_N\| \le C_Nd_S^N$ ($C_N$ is an absolute constant depending on $N$) on $S.$ If 
 \[
 \ dist(g_N, R(K\cap S)) \le c_N d_S^{N},
\]
where 
 \[
 \ c_N = \dfrac{N!}{2^{N+8} \pi C_N^\phi}
 \]
and $C_N^\phi$ is in (2-11). Then 
 \[
 \ \gamma (Int(S)\setminus K) \ge \dfrac{4c_N}{C_N+4c_N}d_S.
 \] 
\end{Lemma}

Proof: 
Let $\phi_1 = \bar\partial^{N-1}\phi(\frac{2z}{d_S}),$ where $\phi$ is in (2-11). Let 
$ f_1 = T_{\phi_1} g_N,$ then
  \[
	\ \| f_1 \| \le 2^{N+3} C_N^\phi C_Nd_S 
	\]
	and
 \[
 \ f_1'(\infty) = \dfrac{1}{\pi} \int \bar\partial g_N \phi_1 dA = (-1)^{N-1} \dfrac{N!d_S^2}{4\pi}.
 \]
 We have the following computation:
 \[
 \ \begin{aligned}
 \ & dist( f_1, R(\mathbb C_\infty \setminus (Int(S)\setminus K))) \\
 \  = & \sup _{\underset{\|\mu \| = 1}{\mu \perp R(\mathbb C_\infty \setminus (Int(S)\setminus K))}} \left |\int f_1 d\mu  \right | \\
\ = & \sup _{\underset{\|\mu \| = 1}{\mu \perp R(\mathbb C_\infty \setminus (Int(S)\setminus K))}} \left |\int g_N dT_{\phi_1}^*\mu \right | \\
\ \le & \sup _{\underset{\|\mu \| = 1}{\mu \perp R(\mathbb C_\infty \setminus (Int(S)\setminus K))}} \| T_{\phi_1}^*\mu \| dist (g_N, R(K\cap S)) \\
\ \le & 8 d_S \|\bar \partial ^N \phi_1 \| dist (g_N, R(K\cap S)) \\
\ \le & 2^{N+3} C_N^\phi c_N d_S,
\ \end{aligned}
 \]
where we use the fact that $T_{\phi_1}^*\mu \perp R(K\cap S)$ for $\mu \perp R(\mathbb C_\infty \setminus (Int(S)\setminus K)).$ 
Let
 \[
 \ f = \dfrac{f_1}{(-1)^{N-1}2^{N+3} C_N^\phi C_Nd_S},
 \]
then $f$ is analytic off $S,$ $\|f\| \le 1, $ $f'(\infty ) = 8\frac{c_N}{C_N}d_S,$ and
 \[ 
 \ dist( f, R(\mathbb C_\infty \setminus (Int(S)\setminus K))) \le \dfrac{c_N}{C_N}.
 \]
Then there is a compact subset $F$ of $Int(S)\setminus K$ and a rational function $r$ with poles in $F$ such that
 \[
 \ \|f - r \|_{\mathbb C_\infty \setminus F} < 2\frac{c_N}{C_N}.
 \]
Hence, $|r(\infty)| < 2\frac{c_N}{C_N}, $ $\|r \|_{\mathbb C_\infty \setminus F} \le 1 + 2\frac{c_N}{C_N}, $ and by the maximum modulus principle,
 \[
 \ |f'(\infty) - r'(\infty) | = \lim_{z\rightarrow \infty} |z(f(z)-r(z)+r(\infty))| \le \max_{|z| = d_S} |z(f(z)-r(z)+r(\infty))| < 4\frac{c_N}{C_N}d_S.
 \]
Therefore,
 \[
 \ \gamma (Int(S)\setminus K) \ge \gamma (F) \ge \dfrac{|r'(\infty) |}{\|r \|_{\mathbb C_\infty \setminus F} + |r(\infty)|}\ge \dfrac{4c_N}{C_N+4c_N}d_S.
 \]

Proof of Theorem \ref{MTheorem1} (2): We assume $\lambda _0 = 0$ and $\delta = 1$ (same argument used in the first paragraph of the proof of Theorem \ref{MTheorem1} (1)). 

We use mathematical induction for $N.$ First we consider the case that $N = 1.$ Let $\epsilon = \dfrac{4c_1}{C_1+4c_1}$ in Lemma \ref{propN}.
 Let $\epsilon _1$ be as in (2-5). Then from  assumption (1-3) and Lemma \ref{LightRoute}, we conclude that Case II of $scheme(Q,\epsilon,m, \gamma_n, \Gamma_n, n\ge m)$ is true. Let $S$ be a heavy $\epsilon$ square and $g_1 = \bar z - c_S,$ then by Lemma \ref{propN}, we must have 
 \[
 \ dist (\bar z, R(S\cap K)) = dist (g_1, R(S\cap K)) \ge c_1 d_S.
 \]
Hence, (2-9) holds for $N=1,$ by Lemma \ref{HBarrier}, we prove (1-5) for $N=1.$

Now we assume that (1-5) holds for $k =1,2,...,N.$ The proof is similar to that of Theorem \ref{MTheorem1} (1).  
Set
 \[
 \ \epsilon = min (c_0, \dfrac{\epsilon _{N}}{2}, \dfrac{\epsilon _{N-1}}{2^2},...,\dfrac{\epsilon _1}{2^{N}},\dfrac{\epsilon _0}{2^{N+1}}),
 \] 
where $c_0 > 0$ will be determined later, and $\epsilon _{N+1}$ as in (2-5). Since $\gamma (B (0, 1) \setminus K) < \epsilon_{N+1}$ (assumption (1-3)), from Lemma \ref{LightRoute}, we conclude that Case II of $scheme(Q,\epsilon,m, \gamma_n, \Gamma_n, n\ge m)$ is true. Let $S$ be a heavy $\epsilon$ square. We assume
 \[
 \ dist \left (\bar z ^{N+1}, HP (1, \bar z, ..., \bar z^N, K \cap S ) \right ) \le\frac{1}{2} d_S^{N+1}, 
 \]
otherwise (2-9) already holds. Without loss of generality, we assume the center of $S$ is zero. There are polynomials $p_1,..., p_N$ and a rational function $r$ with poles off $K\cap S$ such that
 \[
 \ \left \| \bar z^{N+1} + \sum _{j=1}^N p_j \bar z^j + r \right \| _{K\cap S} \le d_S^{N+1}.
 \]
Using the same argument of the paragraph (in the proof of Theorem \ref{MTheorem1}) under (2-12), we get 
 \[
 \ \left \| \bar z^{N+1} + \sum _{j=1}^N p_j \bar z^j \right \| _{K\cap \bar B(0, \frac{d_S}{2^{N+1}})} \le C_{N+1}d_S^{N+1}.
 \]
Let $g_{N+1} = \bar z^{N+1} + \sum _{j=1}^N p_j \bar z^j,$ then $\|g_{N+1}\| \le C_{N+1}d_S^{N+1}$ on $\frac{1}{2^{N+1}} S \subset B(0, \frac{d_S}{2^{N+1}}).$
By Lemma \ref{propN} for $\frac{1}{2^{N+1}} S,$ if we choose
 \[
 \ c_0 = \dfrac{4c_{N+1}}{C_{N+1}+4c_{N+1}},
 \] 
we must have
 \[
 \ dist \left (\bar z ^{N+1}, HP (1, \bar z, ..., \bar z^N, K \cap S ) \right ) = dist \left (g_{N+1}, HP (1, \bar z, ..., \bar z^N, K \cap S ) \right ) \ge c_{N+1} d_S^{N+1}. 
 \]
Therefore, (2-9) holds and the theorem now follows from Lemma \ref{HBarrier}.

Proof of Theorem \ref{BPETheorem}: Assume that $\lambda_0$ is not a peak
point for $R(K),$ then by Curtis's Criterion (see Theorem 4.1 in \cite{gamelin}, p. 204), we get,
 \[
 \ \limsup_{\delta\rightarrow 0} \dfrac{\gamma (B (\lambda _0, \delta) \setminus K)}{\delta} = 0 .
 \]
Then together with the assumptions for $F_1,...,F_N,$ we can choose $\delta > 0$ such that
 \[
 \ F_j(z) = \sum _{i=0}^M g_{ij}(z) \bar z ^i, ~ z\in B (\lambda _0, \delta), 
 \]
and
 \[
 \ \gamma (B (\lambda _0, \delta) \setminus K) < \epsilon \delta  \tag{2-14}
 \]
where $M \ge N,$ $g_{ij} \in A(\bar B (\lambda _0, \delta)),$ $j=1,...,N,$ $i = 0,1,...,M,$
 \[
 \ \epsilon = \min \left(\epsilon _M, \dfrac{\epsilon _{M-1}}{2}, ..., \dfrac{\epsilon _{1}}{2^{M-1}} \right),
 \]
 and $\epsilon _{1},...,\epsilon _{M}$ are in Theorem \ref{MTheorem1}. Set   
 \[
 \ F_0(z) = \sum _{i=0}^M g_{ij}(z) \bar z ^i, ~ z\in B (\lambda _0, \delta), 
 \]
where $g_{00} = 1$ and $g_{i0} = 0$ for $i = 1,...,M.$ Then
 \[
 \ \left({\begin{array}{cccc} F_0 & F_1& ... & F_N \\ 0 & \bar\partial F_1& ... & \bar\partial  F_N \\ ... & ...& ... & ... \\ 0 & \bar\partial ^N F_1& ... & \bar\partial ^N F_N\end{array}}\right)
 \ = 
 \ \left({\begin{array}{ccccc} 1 & \bar z&\bar z^2& ... & \bar z^M \\ 0 & 1& 2\bar z&...& M\bar z^{M-1} \\ 0 & 0& 2&... & M(M-1)\bar z^{M-2} \\ ... &... & ...& ... & ... \\ 0 & 0& 0&... & \frac{M!}{(M-N)!}\bar z^{M-N}\end{array}}\right)
 \ \left({\begin{array}{ccc} g_{00} & ... & g_{0N} \\ g_{10} & ... & g_{1N} \\ ...& ... & ... \\ g_{M0} & ... & g_{MN}\end{array}}\right)
 \]
Since the matrix $[\bar\partial ^iF_j(\lambda_0)]_{1\le i,j\le N}$ is invertible, we may choose above $\delta$ small enough and $i_0=1 <i_1 <i_2 <...<i_N$ so that 
 \[
 \ G(\lambda ) = \left({\begin{array}{cccc} g_{i_00} &g_{i_01} & ... & g_{i_0N}  \\ g_{i_10} &g_{i_11} & ... & g_{i_1N}  \\ ...&...& ... & ... \\g_{i_N0} &g_{i_N1} & ... & g_{i_NN}\end{array}}\right)
 \]
is analytic and invertible on $B (\lambda _0, \delta).$ Moreover, 
 \[
 \ \|G(\lambda )^{-1}\| \le C, ~ \lambda \in B (\lambda _0, \delta).
 \]
So $G(\lambda )$ and $G(\lambda )^{-1}$ are linear and uniform bounded operators from $(\mathbb C^{N+1}, \|.\|_{N+1})$ to $(\mathbb C^{N+1}, \|.\|_{N+1}),$ where $\|x\|_{N+1} = \sum_{i=1}^{N+1} |x_i|.$ For polynomials $p_0,p_1,...,p_N,$ we get
 \[
 \ \begin{aligned}
 \ \sum_{i=0}^M\left |\sum_{j=0}^N p_j(\lambda ) g_{ij}(\lambda ) \right | \ge &\sum_{k=0}^N\left |\sum_{j=0}^N p_j(\lambda ) g_{i_kj}(\lambda ) \right |\\
 \ = & \|(p_0(\lambda ),p_1(\lambda ),...,p_N(\lambda ))G(\lambda )\| \\
 \ \ge &\frac{1}{C} \sum_{j=0}^N |p_j(\lambda )|
 \ \end{aligned} \tag{2-15}
 \] 
on $B (\lambda _0, \delta).$ 

Now let us prove (1). From (2-14) and Theorem \ref{MTheorem1} (1), we have
 \[
 \ |q_M(\lambda )| \le \dfrac{C_M}{\delta ^{M+2}} \left\| \sum _{i=0}^M q_i\bar z^i \right \|_{L^1(K\cap\bar B (\lambda _0, \delta))}
 \]
on $B (\lambda _0, \frac{\delta}{2})$ and $q_0,q_1,...,q_M$ are polynomials. Using (2-14) and Theorem \ref{MTheorem1} (1) again, we have
the following calculation
 \[
 \ \begin{aligned}
 \ |q_{M-1}(\lambda )| \le & \dfrac{C_{M-1}}{(\frac{\delta}{2}) ^{M+1}} \left\| \sum _{i=0}^{M-1} q_i\bar z^i\right \|_{L^1(K\cap\bar B (\lambda _0, \frac{\delta}{2}))} \\
 \ \le & \dfrac{C_{M-1}}{(\frac{\delta}{2}) ^{M+1}} \left ( \left \| \sum _{i=0}^{M} q_i\bar z^i \right\|_{L^1(K\cap\bar B (\lambda _0, \frac{\delta}{2}))} +  \| q_M  \|_{K\cap\bar B (\lambda _0, \frac{\delta}{2})} \|\bar z^M \|_{L^1(K\cap\bar B (\lambda _0, \frac{\delta}{2}))} \right )\\
 \ \le & \dfrac{C_{M-1}}{(\frac{\delta}{2}) ^{M+1}} \left (1 + \dfrac{C_M}{\delta ^{M+2}} \| \bar z^M \|_{L^1(K\cap\bar B (\lambda _0, \delta))} \right )\left\| \sum _{i=0}^{M} q_i\bar z^i \right\|_{L^1(K\cap\bar B (\lambda _0,\delta))}
 \ \end{aligned}
 \]
on $B (\lambda _0, \frac{\delta}{4}).$ Therefore, there is a constant $C_N >0$ such that
\[
 \ \sum _{i=0}^M |q_i(\lambda )| \le C_N \left\| \sum _{i=0}^M q_i\bar z^i \right \|_{L^1(K\cap\bar B (\lambda _0, \delta))}\tag{2-16}
 \]
on $B (\lambda _0, \frac{\delta}{2^{M+1}}).$ From (2-15) and (2-16), we see that
 \[
 \ \sum_{j=0}^N |p_j(\lambda )| \le C \sum_{i=0}^M\left |\sum_{j=0}^N p_j(\lambda ) g_{ij}(\lambda ) \right | \le CC_N \left \|\sum_{j=0}^N p_jF_j \right \|_{L^1(K\cap\bar B (\lambda _0,\delta))} \tag{2-17}
 \]
on $B (\lambda _0, \frac{\delta}{2^{M+1}}).$ So $\lambda_0$ is an analytic bounded point evaluation for $P^t(1, F_1,...,F_N,K\cap\bar B (\lambda _0,\delta)).$

The proof of (2) is the same. (2-16) becomes 
 \[
 \ \sum _{i=1}^M |q_i(\lambda )| \le C_N \left\| \sum _{i=1}^M q_i\bar z^i + r \right \|_{K\cap\bar B (\lambda _0, \delta)}
 \]
where $r$ is a rational function with poles off $K\cap\bar B (\lambda _0, \delta).$ (2-17) becomes 
\[
 \ \sum_{j=1}^N |p_j(\lambda )| \le C \sum_{i=1}^M\left |\sum_{j=1}^N p_j(\lambda ) g_{ij}(\lambda ) \right | \le CC_N \left \|\sum_{j=1}^N p_jF_j +r \right \|_{K\cap\bar B (\lambda _0,\delta)}
 \]
on $B (\lambda _0, \frac{\delta}{2^{M+1}}).$ So $\lambda_0$ is an analytic bounded point evaluation for $HP(F_1,...,F_N,K\cap\bar B (\lambda _0,\delta)).$

\section{Uniform Rational Approximation}

In this section, we will prove Theorem \ref{MTheorem2} and Proposition \ref{Prop}. To prove Theorem \ref{MTheorem2}, we need to prove several Lemmas.

\begin{Lemma}\label{Special}
If $R(K) \ne C(K),$ then $ HP(\bar z,\bar z^2, ..., \bar z^N,K) \ne C(K).$
\end{Lemma}

Proof: Notice, by Lemma \ref{HBarrier} (2-10),  
\[
 \ |p_N(\lambda)| \le \dfrac{C_N}{\delta ^N} \left \|\sum_{k=1}^N p_k \bar z ^k + r \right \|_{K\cap\bar B (\lambda_0, \delta) \setminus B (\lambda_0, \frac{\delta}{2})}, ~\lambda \in  \bar B \left (\lambda_0, \frac{1}{2}\delta \right ).
 \]
Then there exists a finite Borel measure $\mu$ supported on $K\cap\bar B (\lambda_0, \delta) \setminus B (\lambda_0, \frac{\delta}{2})$ such that 
 \[
 \ p_N(\lambda_0) = \int (\sum_{k=1}^N p_k \bar z ^k + r) d\mu ,
 \]
where $r\in Rat(K\cap\bar B (\lambda_0, \delta) \setminus B (\lambda_0, \frac{\delta}{2}))$ and $p_k\in\mathcal P.$ Therefore, the non zero measure $(z-\lambda_0)\mu\perp HP(\bar z,\bar z^2, ..., \bar z^N,K)$ and $HP(\bar z,\bar z^2, ..., \bar z^N,K) \ne C(K).$

\begin{Lemma}\label{OnSquare}
Let $T$ be a closed square and $\gamma(B (c_T, \sqrt{2} d_T) \setminus K) < \epsilon_N \sqrt{2} d_T.$ If $A(K)\subset HP(\bar z,\bar z^2, ..., \bar z^N,K),$ $\phi$ is a smooth function supported inside $T,$ and $f\in A(K),$ then $T_\phi f \in R(K).$ 
\end{Lemma}

Proof: Case 1: Suppose $Int(K\cap T) = \emptyset.$ Then $A(K\cap T) = C(K\cap T)\subset HP(\bar z,\bar z^2, ..., \bar z^N,K\cap T)$ since each smooth function with support in $K\cap T$ belongs to $A(K).$ From Lemma \ref{Special}, we get $C(K\cap T) = R(K\cap T).$ Hence, $T_\phi f\in R(K).$

Case 2: $Int(K\cap T) \ne \emptyset.$  
There are sequences of $\{p_{ij}\}_{1\le i\le N, 1\le j < \infty}\subset \mathcal P$ and $\{r_j\}\subset Rat(K)$ such that 
 \[ 
 \ \lim_{j\rightarrow\infty}(\sum_{i=1}^N p_{ij}\bar z^i + r_j) = T_\phi f
 \]
uniformly on $K$ since $T_\phi A(K) \subset A(K).$ By Theorem \ref{MTheorem1} (2), for each $i$ , the sequence $\{p_{ij}\}_{1\le j < \infty}$ converges to an analytic function $f_i$ uniformly on $T\subset \bar B (c_T, \sqrt{2} d_T/2).$ Hence $\{r_j\}$ converges to $r$ uniformly on $K\cap T$ that is analytic on $Int(K\cap T)$ and $T_\phi f(z) = \sum_{i=1}^N f_i(z)\bar z^i + r(z)$ on $Int(K\cap T).$ This implies $f_i(z) = 0$ on $T$ and $T_\phi f\in R(K\cap T).$ The Lemma is proved. 

\begin{Lemma}\label{ML2}
If $A(K)\subset HP(\bar z,\bar z^2, ..., \bar z^N,K),$ then $A(K) = R(K).$ 
\end{Lemma}

Proof: We use standard Vitushkin approximation scheme (see \cite{gamelin} for example). Let $\{\psi_n,S_n\}$ be a smooth partition of unity, where the length of $S_n$ is $\delta,$ the support of $\psi_n$ is in $2S_n, $ $\|\bar\partial \psi_n \| \le C/\delta,$ $\sum \psi_n = 1,$ and $\cup_{n=1}^\infty S_n = \mathbb C$ with $Int(S_n)\cap Int(S_m) = \emptyset.$  For a function $f\in A(K),$
 \[
 \ f = \sum_{n = 1}^\infty T_{\psi_n}f.
 \]
For a fixed $n,$ let $T = 2S_n,$ if $\gamma(B (c_T, \sqrt{2} d_T) \setminus K) < \epsilon_N \sqrt{2} d_T,$ then, by Lemma \ref{OnSquare}, $h_n = T_{\psi_n}f \in R(K).$ If $\gamma(B (c_T, \sqrt{2} d_T) \setminus K) \ge \epsilon_N \sqrt{2} d_T,$ then $\gamma(Int(2T)\setminus K) \ge \epsilon_N \sqrt{2} d_T.$ Since
\[
 \ \left |\int f \bar\partial \psi_n dA \right | \le C_1 d_T \omega(f,\delta),~ \left |\int (z-c_T)f \bar\partial \psi_n dA \right | \le C_1 d_T ^2 \omega(f,\delta),
 \]
where
 \[
 \ \omega(f,\delta) = \sup_{z,w\in B(c_T, \sqrt{2}\delta)} |f(z) - f(w)|,
 \]
we set 
 \[
 \ \alpha = \dfrac{\int f \bar\partial \psi_n dA }{C_1 d_T \omega(f,\delta)},~ \beta = \dfrac{\int (z-c_T)f \bar\partial \psi_n dA}{ C_1 d_T ^2 \omega(f,\delta)}.
 \] 
Using Lemma \ref{prop1}, we can find a function $g \in R(\mathbb C_\infty \setminus (Int(2T)\setminus K)) \subset R(K)$ satisfying (1) to (5) of Lemma \ref{prop1}. Now let $ h_n = \frac{C_1 \omega(f,\delta)}{\pi}g,$
then $h_n\in R(K),$ $h_n$ is analytic off $2T,$ $\|h_n\| \le C \omega(f,\delta),$ and $h_n - T_{\psi_n}f$ has triple zeros at $\infty.$ So  $\sum _{n=1}^\infty h_n$ goes to $f$ uniformly when $\delta$ tends to zero.  This completes the proof of the lemma.

Proof of Theorem \ref{MTheorem2}: Let $M$ be the largest power of $\bar z$ among all terms of $q_1(z,\bar z),q_2(z,\bar z),...,q_N(z,\bar z),$ then 
 \[
 \ A(K) \subset HP(q_1,...,q_N,K)\subset HP((\bar z,\bar z^2,...,\bar z^M) ,K). 
 \]
By Lemma \ref{ML2}, we conclude $A(K) = R(K).$

Can Theorem \ref{MTheorem2} work for more general functions? For $N=1$ and a smooth function $g$ with $\bar\partial g \ne 0,$ it is proved in \cite{yang2} that $A(K)\subset HP(g,K)$ implies $A(K) = R(K).$ We think this may still hold for $N > 1.$ However, in this section, we provide an example for a single function as stated in Proposition \ref{Prop}.

Let us construct a compact subset $K_0$ of the closed unit square with center at zero and sides parallel to coordinate axes such that $P(K_0) = C(K_0)$ and $Area(K_0) = a, ~ 0 < a < 1.$ In fact, we can construct a planar Cantor set $K_0$ as the following. Given a sequence $\{\lambda _n \}$ with $0 < \lambda _n < \frac{1}{2},$ let $Q_0 = [0, 1] \times [0, 1].$ At the first step we take four closed squares inside $Q_0,$ with side length $\lambda_1,$ with sides parallel to the coordinate axes, and so that each square contains a vertex of $Q_0.$ At the second step we apply the
preceding procedure to each of the four squares obtained in the first step, but now using the proportion factor $\lambda _2.$ In this way, we get 16 squares of side length $\sigma_2 = \lambda_1 \lambda_2.$ Proceeding inductively, at each step we obtain $4^n$ squares $Q^n_j,~ j = 1,2,... 4^n$ with side length $\sigma _n = ¸\lambda_1 \lambda_2...\lambda_n.$ Now let
 \[
 \ L_n = \cup _{j=1}^{4^n} Q^n_j, ~ K_0 = \cap _{n=1}^\infty L_n,
 \]
and
 \[
 \ \lambda _n = \dfrac{1}{2} a^{\frac{1}{2^{n+1}}}.
 \]
Then
 \[
 \ Area(K_0) = \lim_{n\rightarrow \infty} 4^n \sigma_n^2 = a,  
 \]
By construction, $\mathbb C\setminus K_0$ is connected, so 
 \[
 \ P(K_0) = C(K_0).\tag{3-1}
 \]
Since $\lim_{n\rightarrow \infty} \frac{\sigma_n}{\sigma_{n+1}} = 2,$ we can choose $n_0$ so that for $n\ge n_0,$ we have $2 \sigma_{n+1} \le \sigma_n \le 2.1 \sigma_{n+1}.$ Now let $T$ be a square with $d_T \le \sigma_{n_0},$ then we can find an integer $n_1 > n_0$ such that $2 \sigma_{n_1+1} \le\sigma_{n_1} \le d_T \le \sigma_{n_1 - 1}.$ 
Suppose $T\cap K_0 \ne \emptyset,$ there exists $Q^{n_1+1}_{j_1}$ with $T\cap Q^{n_1+1}_{j_1} \ne \emptyset. $ Therefore, $Q^{n_1+1}_{j_1} \subset 2T$ and
 \[
 \ Area ((2T)\cap K_0) \ge Area (Q^{n_1+1}_{j_1}\cap K_0) = \lim_{n\rightarrow \infty} 4^n \sigma_{n_1+1+n}^2 \ge a\sigma_{n_1+1}^2 \ge \dfrac{a}{(2(2.1)^2)^2} Area ((2T).
 \]
In this case,
 \[
 \ \gamma ((2T)\cap K_0)) \ge \dfrac{a}{(4(2.1)^2)\sqrt{\pi}} (=c_0) d_T. \tag{3-2}
 \]
Now we will construct a sequence of disjoint small open disks $\{B_k(z_k, r_k) \}_{k=1}^\infty$ within $G = B(0, 1)\setminus K_0$ satisfying the following conditions:

(1) Each point in $K_0$ is a limit of a subsequence of the disks;

(2) No point in $G$ is a limit of a subsequence of the disks;

(3) $\sum_{k=1}^\infty r_k < \infty.$

Let $\{x_k\} \subset K_0$ be a subset that is dense in $K_0.$ We begin with a point $y_{11}\in G,$ choose $0 < r_{11} < \frac{1}{2}$ with $B(y_{11},r_{11})\subset G.$ Let $d_{11} = dist (y_{11}, K_0) - r_{11}.$ We finish the level 1 construction. For level 2, choose $y_{21}\in G$ so that $dist(y_{21}, x_1) < \frac{d_{11}}{2}$ and $ 0 < r_{21} <  \min (\frac{d_{11}}{4}, \frac{1}{2^2})$ with $B(y_{21},r_{21})\subset G.$ Define $d_{21} = dist (y_{21}, K_0) - r_{21}.$ Choose $y_{22}\in G$ so that $dist(y_{22}, x_2) < \frac{d_{21}}{2}$ and $ 0 < r_{22} <  \min (\frac{d_{21}}{4}, \frac{1}{2^2})$ with $B(y_{22},r_{22})\subset G.$ Define $d_{22} = dist (y_{22}, K_0) - r_{22}.$ We continue this process and get $\{y_{ij}\},$ $\{r_{ij}\},$ and $\{d_{ij}\}$ satisfying:
 \[
 \ dist (y_{ij}, K_0) = d_{ij} + r_{ij},~ \sum r_{ij} < \infty;
 \] 
and
 \[
 \  d_{i,1} < \dfrac{d_{i-1,i-1}}{2}, ~ d_{ij} <  \dfrac{d_{i,j-1}}{2},
 \]
where $1<j\le i.$ Let $ \{z_k\} = \{y_{ij}\}$ and $\{r_k\} = \{r_{ij}\}.$ Clearly the conditions (1)-(3) are met. Set
 \[
 \ K = \bar B(0,1) \setminus \left(\bigcup_{k=1}^\infty B(z_k,r_k) \right). \tag{3.3}
 \]
It is easy to verify that the inner boundary $\partial_i K$ equals $  K_0.$ For this $K,$ we can prove Proposition \ref{Prop}.

Proof of Proposition \ref{Prop}: Let $\phi = \mathcal C(dA|_{K_0}),$ then $\phi\in A(K).$ Define 
 \[
 \ \mu = dz |_{\partial B (0,1)} - \sum_{k=1}^\infty dz |_{\partial B (z_k,r_k)},
 \]
 then $\mu$ is a finite Borel measure, $\mu \perp R(K),$ and 
 \[
 \ \int \phi d\mu = \int_{K_0}\int\dfrac{1}{\lambda - z} d\mu (z) dA(\lambda) = - 2 \pi i Area (K_0) \ne 0.
 \] 
Hence, $R(K) \ne A(K).$ Since for a polynomial $p,$
 \[
 \  \mathcal C(pdA|_{K_0}) - p\phi = \int_{K_0} \dfrac{p(w) - p(z)}{w - z} dA(w) \in R(K),
 \]
we see that $\mathcal C (pdA|_{K_0})\in HP(\phi, K).$ Using (3-1), we conclude that 
 \[
 \ \mathcal C(\chi_E dA) \in HP(\phi,K), ~ E \subset K_0, \tag{3-3}
 \]
 where $\chi_E$ is the characteristic function of $E.$ Let $T$ and $Q^{n_1 +1}_{j_1}$ be squares as above. There are four squares 
 \[
 \ \cup_{i=1}^4Q^{n_1 +2}_{k_i}\subset Q^{n_1 +1}_{j_1}.
 \]
Choose two of them, say $Q^{n_1 +2}_{k_1}$ and $Q^{n_1 +2}_{k_2},$ with the same y coordinates of the centers. Set $E_1 = Q^{n_1 +2}_{k_1} \cap K_0$ and $E_2 = Q^{n_1 +2}_{k_2} \cap K_0.$ Let $\phi_1 (z) = - \chi_{E_1}$ and $\phi_2 (z) = - \chi_{E_2}.$ Set $f_j = \mathcal C(\phi_j dA)$ for $j = 1,2.$ Then from (3-3), we have $f_j\in HP(\phi, K)$ and 
 \[
 \ f_j'(\infty ) = Area(E_j) (=AE)
 \]
since $Area(E_1) = Area(E_2),$ and
 \[
 \ \beta (f_j, c_T) = \int (z - c_{E_j}) \chi_E (z)dA + (c_{E_j} - c_T )Area(E_j) = (c_{E_j} - c_T )AE.
 \]
We Set 
 \[
 \ g_1 = \dfrac{(c_{E_2} - c_T )f_1 - (c_{E_1} - c_T)f_2}{(c_{E_2} - c_{E_1} )AE}d_T,
 \]
 and
 \[
 \ g_2 = \dfrac{f_2 - f_1}{(c_{E_2} - c_{E_1} )AE} d_T^2 .
 \]
 Notice that $|c_{E_2} - c_T| < 2d_T, $ $|c_{E_1} - c_T| < 2d_T,$ and $c_{E_2} - c_{E_1}$ is comparable with $d_T.$ Using the arguments before (3-2), we see that $\|g_j\| \le C$ ( absolute constant), $g_j\in HP(\phi, K),$ and 
 \[
 \ g_1'(\infty ) = d_T, ~ g_2'(\infty ) = 0, ~ \beta (g_1, c_T) = 0, ~ \beta (g_2, c_T) = d_T^2.
 \]
Now we use standard Vitushkin approximation scheme (see \cite{gamelin} for example). Let $\{\psi_n,S_n\}$ be a smooth partition of unity, where the length of $S_n$ is $\delta,$ the support of $\psi_n$ is in $2S_n, $ $\|\bar\partial \psi_n \| \le C/\delta,$ $\sum \psi_n = 1,$ and $\cup_{n=1}^\infty S_n = \mathbb C$ with $Int(S_n)\cap Int(S_m) = \emptyset.$ We assume $\delta$ is less than $\frac{\sigma_{n_0}}{2}.$ For a function $f\in A(K),$
 \[
 \ f = \sum_{n = 1}^\infty T_{\psi_n}f.
 \]
For a fixed $n,$ if $(2S_n) \cap K_0 = \emptyset,$ then $h_n = T_{\psi_n}f \in R(K).$ If $(2S_n) \cap K_0 \ne \emptyset,$ then let $T=2S_n$ and
\[
 \ h_n = \dfrac{\int f \bar\partial \psi_n dA}{\pi d_T}g_1 + \dfrac{\int (z-c_T)f \bar\partial \psi_n dA}{\pi d_T^2}g_2, 
 \]
then $h_n\in HP(\phi, K),$ $h_n$ is analytic off $2T,$ $\|h_n\| \le C \omega(f,\delta),$ and $h_n - T_{\psi_n}f$ has triple zeros at $\infty.$ So $\sum _{n=1}^\infty h_n$ goes to $f$ uniformly when $\delta$ tends to zero.  This completes the proof of the proposition.

\bibliography{Bibliography}
\end{document}